\documentclass[12pt,a4]{article}
\usepackage{setspace}
\usepackage{amssymb}
\usepackage{amsmath,amssymb}
\usepackage[T1]{fontenc}
\usepackage[bf,small,tableposition=top]{caption}
\usepackage{subfig}
\usepackage{amsthm}
\usepackage{lineno}
\usepackage[figuresright]{rotating}
\usepackage{enumerate}
\usepackage{authblk}
\usepackage[T1]{fontenc}
\usepackage[utf8]{inputenc}
\usepackage{graphics}
\usepackage{graphicx}
\usepackage{multirow}
\usepackage{lscape}
\usepackage[dvipsnames]{xcolor}
\usepackage{mathrsfs}
\usepackage{hyperref}
\usepackage{xcolor}
\usepackage{textcomp}
\usepackage{listings}
\usepackage[dvipsnames]{xcolor}
\begin{document}
	\title{A New Mixed Generalized Negative Binomial Distribution}
	\author{Anwar Hassan$^{1}$ , Ishfaq S. Ahmad$^{2}$ and Peer Bilal Ahmad$^{3\dag}$\\
	
$^{1,2}$  P. G.Department of Statistics, University of Kashmir, Srinagar, India \\
$^\dag$ Corresponding author- peerbilal@yahoo.co.in \\
$^{3}$ Department of Mathematical Sciences, IUST, Srinagar, India
}

\maketitle

\begin{abstract}
	In this article, a generalized version of Negative binomial- beta exponential distribution with five parameters have been introduced. Some interesting submodels have been derived from it. A comprehensive mathematical treatment of  proposed distribution is being provided. Various  expressions like that of  moment generating function, moments are derived. The model parameters are estimated by the maximum likelihood method. Finally, the application of proposed distribution is carried out on one sample of automobile insurance data.
\end{abstract}


\section{Introduction}

	Count data often appear in disciplines such as insurance, health, environment, ecology, actuarial science, biology etc. The examples of count datasets, includes, the number of times a patient is referred to  hospital in a calendar year, the number of days a patient stays in the hospital, the number of accidents, the number of children a married couple has, the number of insurance claims, the number of kinds of species in ecology and many more. The typical features possessed by such datasets are over dispersion, zero vertex unimodality and positively skewed nature. In many cases, the simple Poisson distribution model is not appropriate because it imposes the restriction that the conditional mean of each count variable must equal the conditional variance.
	As Min and Czado (2010) point out, the Poisson model  is too simple to capture complex	structures of count data such as over dispersion and zero vertex modalitiy, that is , presence of higher percentage of zero values. In consequence, many attempts have been made to develop such models that are less restrictive than Poisson, and are based on other distributions, have been presented in the statistical literature, including the negative binomial, generalized Poisson and generalized negative binomial (see Cameron and Trivedi (1998) and Famoye (1995), among others). Also various methods have been employed to develop new class of discrete distributions like  mixed Poisson method (see Karlis and Xekalaki, 2005), mixed Negative Binomial (Gomez et al., 2008), discretization of continuous family of distribution and discrete analogues of continuous distribution. \par 

	 Jain and Consul (1971) proposed a generalized  negative binomial distribution   i.e., $ X \sim \mathcal{GNB}(m,\theta, \beta )$ with probability mass function (PMF)  
	\begin{equation}
	p(x)=\frac{m}{m+\beta x}\binom{m+\beta x}{x}(1-\theta)^{x} \theta^{m+\beta x-x},\quad x=0,1,\dots
	\end{equation}
 
where $m>0, 0  \leq \theta <1, \theta \beta <1$ and $\beta \geq 1$. Let $\theta=e^{-\lambda}$ follows  Beta Exponential distribution proposed by Nadarajah and Kotz (2006) denoted as $ \theta \sim \mathcal{BE}(a,b,c)$ with pdf given by 
	\begin{equation} \label{be}
g(x)=\frac{c}{B(a,b)}e^{-bcx}\left[1-e^{-cx}\right]^{a-1},\quad x>0 
\end{equation}
where $ a,b,c>0$ and   $B(r,s)=\frac{\Gamma(r)\Gamma(s)}{\Gamma(r+s)}$; $r,s>0$. \par 
The moment generating function (mgf) corresponding to (\ref{be}) is given by 
\begin{equation}
M_{x}(t)=\frac{B(b-\frac{t}{c},a)}{B(a,b)}
\end{equation}

The aim of our work is to introduce generalized version of Negative Binomial-Beta Exponential distribution (Pudprommart et al.,2012) and discuss its various special cases which are actually some well known models. The factorial moments corresponding to those distributions are also discussed. Finally a data set has been used to check the efficacy of our proposed model.

\section{A New Mixture distribution}
The PMF of proposed distribution is defined by following stochastic representation
\begin{equation}  \label{sr}
\begin{aligned}
X|\lambda\sim  &\mathcal{GNB}(m,\beta, \theta=e^{-\lambda})\\
\lambda|a,b,c  \sim &\mathcal{BE}(a,b,c)
\end{aligned} 
\end{equation}
\textit{where,} $\mathcal{GNB}(m,\beta, \theta=e^{-\lambda})$ \textit{is the generalized negative binomial distribution and} $\mathcal{BE}(a,b,c)$ \textit{Beta Exponential  distribution} and is obtained in Theorem 1. \\

\noindent \textbf{Theorem 1.} \textit{Let $X\sim \mathcal{GNB-BE}(m,\beta,\theta=e^{-\lambda},a,b,c)$ be a Generalized negative binomial-Beta exponential  distribution as defined in (\ref{sr}) then pmf is given by}
\begin{equation} \label{pmf} 
p(x)=\frac{m}{m+\beta x}\binom{m+\beta x}{x}\sum_{j=0}^{x}\binom{x}{j}(-1)^{j}\frac{B(b+\frac{j+m+\beta x-x}{c},a)}{B(a,b)}, \qquad  x=0,1,\cdots,
\end{equation}
\textit{with  $m>0, \beta\geq 1$ and $a,b,c>0$.}

\noindent \textit{Proof:} If $X|\lambda \sim \mathcal{GNB}(m,\beta, \theta=e^{-\lambda}) $ and $\lambda \sim \mathcal{BE}(a,b,c) $ , then pmf of $X$ can be obtained by 
\begin{equation}
h(x)=\int_{0}^{\infty}f(x|\lambda)g(\lambda;a,b,c)d\lambda
\end{equation}
where 
\begin{equation}
\begin{aligned}
f(x|\lambda)=&\frac{m}{m+\beta x}\binom{m+\beta x}{x}e^{-\lambda(m+\beta x-x)}(1-e^{-\lambda})^x \\
=&\frac{m}{m+\beta x}\binom{m+\beta x}{x}\sum_{j=0}^{x}\binom{x}{j}(-1)^{j}e^{-\lambda(j+m+\beta x-x)}
\end{aligned}
\end{equation}
Putting Equation $(7)$ in Equation $(6)$, we find
\begin{equation}
\begin{aligned}
h(x)=&\frac{m}{m+\beta x}\binom{m+\beta x}{x}\sum_{j=0}^{x}\binom{x}{j}(-1)^{j}\int_{0}^{\infty} e^{-\lambda(j+m+\beta x-x)}g(\lambda;a,b,c)d\lambda\\
=&\frac{m}{m+\beta x}\binom{m+\beta x}{x}\sum_{j=0}^{x}\binom{x}{j}(-1)^{j}M_{\lambda}\left(-(j+m+\beta x-x)\right)
\end{aligned}
\end{equation}
Putting the mgf of $\mathcal{BE}(a,b,c)$ in Equation $(3)$ into Equation $(8)$, we get the pmf of $\mathcal{GNB-BE}(m,\beta,a,b,c)$ as 
\begin{equation*} 
p(x)=\frac{m}{m+\beta x}\binom{m+\beta x}{x}\sum_{j=0}^{x}\binom{x}{j}(-1)^{j}\frac{B(b+\frac{j+m+\beta x-x}{c},a)}{B(a,b)}, 
\end{equation*}
which proves the theorem. \hfill $\blacksquare$ 
\begin{figure}[htp]
	\begin{center}
		\includegraphics[scale=0.9]{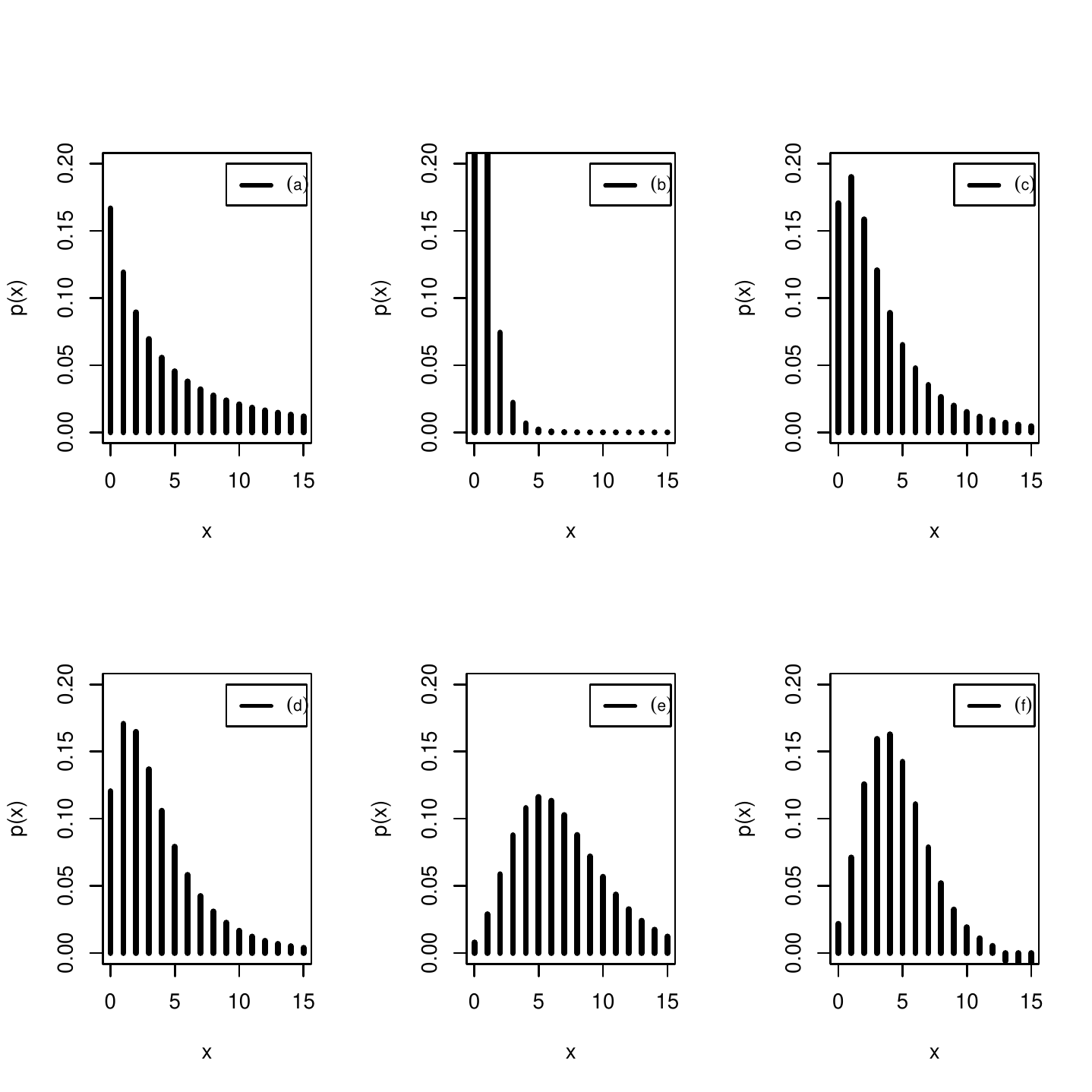}
		\captionof{figure}{PMF plot of $\mathcal{GNB-BE}(m,\beta,a,b,c)$ distribution for different valaues of parameter: (a) $\beta=1,m=1,a=5,b=1,c=1$, (b)$\beta=1,m=3,a=4,b=4,c=5$, (c) $\beta=1, m=3,a=5,b=2,c=2$, (d) $\beta=1,m=5,a=5,b=10,c=0.8$ (e) $\beta=1,m=20,a=15,b=5,c=5$ and (f) $\beta=1, m=30,a=20,b=20,c=5$.}
	\end{center}
\end{figure}
\newpage
\subsection{Relation with other existing distributions}
\begin{enumerate}
	\item For $\beta=1$, it can be easily verified that $\mathcal{GNB-BE}(m,\beta,a,b,c)$  reduces to Negative Binomial-Beta Exponential distribution $\mathcal{NB-BE}(m,a,b,c)$ proposed by Pudprommarat et al.(2012) whose pmf is given by 
	\begin{equation}
	h_{1}(x)=\binom{m+ x-1}{x}\sum_{j=0}^{x}\binom{x}{j}(-1)^{j}\frac{B(b+\frac{m+j}{c},a)}{B(a,b)}
	\end{equation}
	\item  \textit{If $\beta=1$, $c=1$, then $\mathcal{GNB-BE}(m,\beta,a,b,c)$ reduces to generalized Waring distribution with pmf obtained as:}
	\begin{equation}
	h_{2}(x)=\frac{\Gamma(a+b)\Gamma(m+b)a_{(x)}m_{(x)}}{\Gamma(b)\Gamma(m+a+b)(m+a+b)_{(x)}}\frac{1}{x!}; \quad x=0,1,\cdots for \quad m,a,b>0
	\end{equation}
\textit{where $m_{(s)}=\frac{\Gamma(m+s)}{\Gamma(m)}$; $m,s>0$}.

\noindent \textit{Proof:} If $X|\lambda \sim  \mathcal{GNB}(m,\beta=1, \theta=e^{-\lambda})$ and  $\lambda \sim  \mathcal{BE}(a,b,c=1)$, then the pmf of $X$ is obtained as 
\begin{equation}
\begin{aligned}
h_{2}(x)=&\binom{m+ x-1}{x}\sum_{j=0}^{x}\binom{x}{j}(-1)^{j}\frac{B(b+m+j,a)}{B(a,b)}\\
=& \binom{m+ x-1}{x}\left(\frac{\Gamma(a)}{B(a,b)}\right)\sum_{j=0}^{x}\binom{x}{j}(-1)^{j}\frac{\Gamma(j+m+b)}{\Gamma(j+m+a+b)}
\end{aligned}
\end{equation}
Using the result of Gardshteyn and Ryzhikh (2007), the sum of binomial terms in Equation $(11)$ is of the form
\begin{equation}
\sum_{j=0}^{x}\binom{x}{j}(-1)^{j}\frac{\Gamma(j+b)}{\Gamma(j+a)}=
\frac{B(x+a-b,b)}{\Gamma(a-b)}
\end{equation}
Therefore, $h_{2}(x)$ can be written as 
\begin{equation*}
\begin{aligned}
h_{2}(x)=&\binom{m+ x-1}{x}\left(\frac{\Gamma(a)}{B(a,b)}\right)\frac{B(x+a,m+b)}{\Gamma(a)}\\
=&\frac{\Gamma(a+b)\Gamma(m+b)}{\Gamma(b)\Gamma(m+a+b)}\frac{\Gamma(a+x)}{\Gamma(m)}\frac{\Gamma(m+x)}{\Gamma(m)}\frac{\Gamma(m+a+b)}{\Gamma(m+a+b+x)}\frac{1}{x!}\\
=&\frac{\Gamma(a+b)\Gamma(m+b)}{\Gamma(b)\Gamma(m+a+b)}\frac{a_{(x)}m_{(x)}}{(m+a+b)_{(x)}}\frac{1}{x!}
\end{aligned}
\end{equation*}
Hence proved. \hfill $\blacksquare$
\item  \textit{If $\beta=1$, $a=1$, $b=k-m$ and $c=1$, then  $\mathcal{GNB-BE}(m,\beta,a,b,c)$ reduces to Waring distribution with pmf given by }
\begin{equation}
h_{3}(x)=\frac{(k-m)\Gamma(m+x)\Gamma(k)}{\Gamma(m)\Gamma(k+x+1)} ; \quad x=0,1,\cdots. \quad  m>0 \quad  \& \quad  k>m
\end{equation}
\noindent \textit{Proof:} Put a=1 $\&$ b=k-m in Equation $(12)$, we get the pmf of $X$ as 
\begin{equation*}
\begin{aligned}
h_{3}(x)=&\frac{\Gamma(k-m+1)\Gamma(K)1_{(x)}r_{(x)}}{\Gamma(k-m)\Gamma(k+1)(m+1)_{(x)}}\frac{1}{x!} \\
=&\frac{\Gamma(k-m+1)\Gamma(K) \Gamma(x+1)\Gamma(m+x)}{\Gamma(k-m)\Gamma(k+1)\Gamma(1)\Gamma(m)}\frac{\Gamma(k+1)}{\Gamma(k+1+x)}\frac{1}{\Gamma(x+1)}\\
=&\frac{(k-m)\Gamma(m+x)\Gamma(k)}{\Gamma(m)\Gamma(k+x+1)}
\end{aligned}
\end{equation*}
Hence proved. \hfill $\blacksquare$
\item \textit{ If $m=1$, $\beta =1$, $a=1$, and $c=1$, then $\mathcal{GNB-BE}(m,\beta,a,b,c)$ reduces to Yule distribution with pmf given by }
\begin{equation}
\begin{aligned}
h_{4}(x)=\frac{bx!}{(b+1)_{(x+1)}} \quad; x=0,1,\cdots, b>0
\end{aligned}
\end{equation}
\noindent \textit{Proof:}  Put $m=1$ and $a=1$ in Equation $(12)$ , the pmf of the $X$ as  
\begin{equation*}
\begin{aligned}
h_{4}(x)=&\frac{B(x+1,b+1)}{B(1,b)}\\
=&\frac{\Gamma(x+1)\Gamma(b+1)\Gamma(b+1)}{\Gamma(x+b+2)\Gamma(b)} \\
=&\frac{bx!}{(b+1)_{(x+1)}}
\end{aligned}
\end{equation*}
Hence proved. \hfill $\blacksquare$
\item For $\beta=0$ and $m \in \mathcal{N}$, then $\mathcal{GNB-BE}(m,\beta,a,b,c)$ reduces to mixture of Binomial distribution with Beta Exponential distribution whose pmf is given by 
\begin{equation}
h_{5}(x)=\binom{m}{x}\sum_{j=0}^{x}\binom{x}{j}(-1)^{j}\frac{B(b+\frac{j+m-x}{c},a)}{B(a,b)}
\end{equation}
\end{enumerate}
\section{Factorial moments and ordinary (crude) moments of a mixture of Generalized Negative Binomial distribution with Beta Exponential distribution}
Note that the factorial moment of  $\mathcal{GNB}(m,\beta,\theta)$ of order $k$ is given by 
\begin{equation*}
\mu _{[k]}(x)=E\left[x(x-1)\cdots(x-k+1)\right]=\frac{\Gamma(m+k)}{\Gamma(m)}\frac{(1-\theta \beta)^k}{\theta^k}
\end{equation*}
\noindent \textbf{Theorem 2.} \textit{If  $X\sim \mathcal{GNB-BE}(m,\beta,\theta=e^{-\lambda},a,b,c)$, then factorial moment of order $k$ is given by}
\begin{equation}
\mu _{[k]}(x)=\frac{\Gamma(m+k)}{\Gamma(m)}\sum_{j=0}^{k}\binom{k}{j}(-\beta)^{j}\frac{B(b-\frac{k-j}{c},a)}{B(a,b)}
\end{equation}
 \textit{Proof:} If  $X|\lambda\sim \mathcal{GNB}(m,\beta,\theta=e^{-\lambda})$ and  $\lambda\sim \mathcal{BE}(a,b,c)$, then factorial moment of order $k$ can be find out by using concept of conditional moments as 
 \begin{equation*}
 \mu _{[k]}(x)=E_{\lambda}\left[\mu _{[k]}(x|\lambda)\right]
 \end{equation*}
 Using the factorial moment of order $k$ of $\mathcal{GNB}(m,\beta,\theta)$ , $\mu _{[k]}(x)$ becomes 
 \begin{equation*}
 \mu _{[k]}(x)=E_{\lambda}\left[\frac{\Gamma(m+k)}{\Gamma(m)}(e^\lambda -\beta)^k \right]=\frac{\Gamma(m+k)}{\Gamma(m)}E_{\lambda}(e^\lambda -\beta)^k
 \end{equation*}
 Through the  binomial expansion of  $(e^\lambda -\beta)^k=\sum_{j=0}^{k}\binom{k}{j}(-\beta)^{j}e^{\lambda(k-j)}$, $\mu _{[k]}(x)$ can be written as 
 \begin{equation*}
 \begin{aligned}
 \mu _{[k]}(x)=&\frac{\Gamma(m+k)}{\Gamma(m)}\sum_{j=0}^{k}\binom{k}{j}(-\beta)^{j}E_{\lambda}(e^{\lambda(k-j)}\\
 =&\frac{\Gamma(m+k)}{\Gamma(m)}\sum_{j=0}^{k}\binom{k}{j}(-\beta)^{j}M_{\lambda}(k-j)
 \end{aligned}
 \end{equation*}
 From the mgf of $\mathcal{BE}(a,b,c)$ in Equation $(3)$ with $t=k-j$, we get finally factorial moment of order $k$ as:
 \begin{equation*}
 \mu _{[k]}(x)=\frac{\Gamma(m+k)}{\Gamma(m)}\sum_{j=0}^{k}\binom{k}{j}(-\beta)^{j}\frac{B(b-\frac{k-j}{c},a)}{B(a,b)}
 \end{equation*}
 which proves the theorem. \hfill $\blacksquare$

\noindent \textbf{Corollary 1.} \textit{If  $\beta=1$, then factorial moment of order $k$ of $\mathcal{GNB-BE}(m,\beta,\theta=e^{-\lambda},a,b,c)$  reduces to factorial moment of Negative Binomial-Beta Exponential distribution, which is computed by substituting $\beta=1$ in Equation $(16)$ and after simplification:}
\begin{equation}
\mu _{[k]}(x)=\frac{\Gamma(m+k)}{\Gamma(m)}\sum_{j=0}^{k}\binom{k}{j}(-1)^{j}\frac{B(b-\frac{k-j}{c},a)}{B(a,b)}
\end{equation}
\noindent \textbf{Corollary 2.} \textit{If  $\beta=1$ and $c=1$, then factorial moment of order $k$ of $\mathcal{GNB-BE}(m,\beta,\theta=e^{-\lambda},a,b,c)$  reduces to:}
\begin{equation}
\mu _{[k]}(x)=\frac{m_{(k)}a_{(k)}}{(b-1)(b-2)\cdots(b-k)}, \quad k=1,2,\cdots \quad for \quad m,a,b>0  \; \& \; b>k
\end{equation}
which is the actually the factorial moment of order $k$ of generalized Waring distribution.\\
\noindent \textit{Proof:}  Putting $c=1$ in Equation $(17)$, we have 
\begin{equation}
\begin{aligned}
\mu _{[k]}(x)=&\frac{\Gamma(m+k)}{\Gamma(m)}\sum_{j=0}^{k}\binom{k}{j}(-1)^{j}\frac{B(b-k+j,a)}{B(a,b)} \\
=&\frac{\Gamma(m+k)}{\Gamma(m)}\sum_{j=0}^{k}\binom{k}{j}(-1)^{j}   \frac{\Gamma(a)}{B(a,b)} \sum_{j=0}^{k}\binom{k}{j}(-1)^{j}  \left( \frac{\Gamma(j+b-k)}{\Gamma(j+b-k+a)}  \right)
\end{aligned}
\end{equation}
By using expansion given in Equation $(12)$, Equation $(19)$ reduces to 
\begin{equation*}
\begin{aligned}
\mu _{[k]}(x)=&\frac{\Gamma(m+k)}{\Gamma(m)}\frac{\Gamma(a)}{B(a,b)}\frac{B(a+k,b-k)}{\Gamma(a)} \\
=&\frac{\Gamma(m+k)}{\Gamma(m)} \frac{\Gamma(a+k)}{\Gamma(a)}\frac{\Gamma(b-k)}{\Gamma(b)} \\
=&\frac{m_{(k)}a_{(k)}}{(b-1)(b-2)\cdots(b-k)}
\end{aligned}
\end{equation*}
\noindent \textbf{Corollary 3.} \textit{If  $\beta=1$, $a=1$, $b=n-m$ and $c=1$, then factorial moment of order $k$ of $\mathcal{GNB-BE}(m,\beta,\theta=e^{-\lambda},a,b,c)$  reduces to the factorial moment of order $k$ of Waring distribution which is given by:}
\begin{equation}
\mu _{[k]}(x)=\frac{m_{(k)}k!}{(n-m-1)(n-m-2)\cdots(n-m-k)}, \; k=1,2,\cdots \; for \; m>0 \; \&\; n-m>k
\end{equation}
\noindent \textit{Proof:} Substituting $a=1$ $\&$ b=n-m in Equation  $(18)$, we have:
\begin{equation*}
\begin{aligned}
\mu _{[k]}(x)=&\frac{m_{(k)}1_{(k)}}{(n-m-1)(n-m-2)\cdots(n-m-k)}\\
=&\frac{m_{(k)}k!}{(n-m-1)(n-m-2)\cdots(n-m-k)}
\end{aligned}
\end{equation*}
\noindent \textbf{Corollary 4.} \textit{If  $\beta=1$, $m=1$, $a=1$ and $c=1$, then factorial moment of order $k$ of $\mathcal{GNB-BE}(m,\beta,\theta=e^{-\lambda},a,b,c)$  reduces to the factorial moment of order $k$ of Yule distribution which is given by:}
\begin{equation}
\mu _{[k]}(x)=\frac{(k!)^2}{(b-1)(b-2)\cdots(b-k)} \quad k=1,2,\cdots  for \; b>0 \; and \; b>k
\end{equation}
The first four moments now can be easily deduced from factorial moments of $\mathcal{GNB-BE}(m,\beta,\theta=e^{-\lambda},a,b,c)$ from  Equation $(16)$ to get mean, variance,  Index of dispersion etc. \par

The mean and variance  is obtained as:
\begin{equation}
E(X)=\frac{mB(b-\frac{1}{c},a)-m\beta B(a,b)}{B(a,b)}
\end{equation}
\begin{equation*}
E(X^2)=\frac{1}{B(a,b)}\left[m(m+1)B(b-\frac{2}{c},a)-(2m^2 +m)\beta B(b-\frac{1}{c},a)+m^2 \beta ^2 B(a,b)\right]
\end{equation*}
\begin{equation}
\begin{split}
V(X)=&\frac{1}{B(a,b)^2}\left[m(m+1)B(b-\frac{2}{c},a)B(a,b)+mB(b-\frac{1}{c},a)B(a,b)\right] \\
 & -\frac{1}{B(a,b)^2} \left[m^2 B(b-\frac{1}{c},a)^2 -2 \beta m B(b-\frac{1}{c},a) B(a,b) \right]
\end{split}
\end{equation}

\section{Estimation}
	In this Section, we will discuss one of the popular method of estimation namely Maximum Likelihood Estimation (MLE) for the estimation of the parameters of $\mathcal{GNB-BE}(m,\beta,a,b,c)$ distribution. Suppose $\underline{\mathbf{x}}=\lbrace x_{1}, x_{2},\cdots,x_{n} \rbrace$  be a random sample of size $n$ from the  $\mathcal{GNB-BE}(m,\beta,a,b,c)$ distribution with pmf (\ref{pmf}). The likelihood function is given by 
		\begin{equation} \label{likeli} 
	L(\beta,m,a,b,c|\underline{\mathbf{x}})=\prod\limits_{i=1}^{n}\frac{m}{m+\beta x_{i}}\binom{m+\beta x_{i}}{x_{i}}\sum_{j=0}^{x_i}\binom{x_i}{j}(-\beta)^j\frac{B(b+\frac{j+m+\beta x_i-x_i}{c},a)}{B(a,b)}
	\end{equation}
	The log-likelihood function corresponding  to (\ref{likeli}) is obtained as 
	\begin{equation}
\begin{split}
\log L(\beta,m,a,b,c|\underline{\mathbf{x}})=&\sum_{i=1}^{n} \log\left(\frac{m}{m+\beta x_{i}}\binom{m+\beta x_{i}}{x_{i}}\right)+ \\ &+\sum_{i=1}^{n}\left[\log\sum_{j=1}^{x_{i}}\binom{x_{i}}{j}(-\beta)^j  \frac{\Gamma(b+\frac{j+m+\beta x_{i}-x_{i}}{c})\Gamma(a+b)}{\Gamma(a+b+\frac{j+m+\beta x_{i}-x_{i}}{c})\Gamma(b)}\right] 
\end{split}
\end{equation}
	The ML Estimates $\hat{\beta}$ of $\beta$, $\hat{m}$ of $m$, $\hat{a}$ of $a$, $\hat{b}$ of $b$ and $\hat{c}$ of $c$, respectively, can be obtained by solving equations 
\begin{equation*}
\frac{\partial \log L}{\partial \beta}=0, \quad \frac{\partial \log L}{\partial m}=0, \quad \frac{\partial \log L}{\partial a}=0 , \quad \frac{\partial \log L}{\partial b}=0 \quad  \text{and} \quad  	\frac{\partial \log L}{\partial c}=0.
\end{equation*}
where 
\begin{equation}
\begin{split}
\frac{\partial \log L}{\partial \beta}=&\sum_{i=1}^{n}\frac{x_i (b (\beta -1) (c+x_i)+(\beta -1) m+x_i-1)}{(\beta -1) (b (\beta -1) (c+x_i)+(\beta -1) m+x_i)}\\
&+\sum_{i=1}^{n}\frac{x_i (-(m+\beta  x_i) \psi ^{(0)}(m+x_i (\beta -1)+1)+(m+\beta  x_i) \psi ^{(0)}(m+x_i \beta +1)-1)}{m+\beta  x_i} ,
\end{split}
\end{equation}
\begin{equation}
\begin{split}
\frac{\partial \log L}{\partial m}=&\sum_{i=1}^{n}\frac{a (\beta -1) c-\beta  x_i}{(a c+b (c+x_i)+m-x_i) (b (\beta -1) (c+x_i)+(\beta -1) m+x_i)}\\
&+\sum_{i=1}^{n}\frac{-m (m+\beta  x_i) \psi ^{(0)}(m+x_i (\beta -1)+1)+m (m+\beta  x_i) \psi ^{(0)}(m+x \beta +1)+\beta  x_i}{m (m+\beta  x_i)} ,
\end{split}
\end{equation}
\begin{equation}
\begin{split}
\frac{\partial \log L}{\partial a}=&\sum_{i=1}^{n}\frac{(m+(-1+b) x_i)}{((a+b) (a c+m-x_i+b (c+x_i)))} ,
\end{split}
\end{equation}
\begin{equation}
\begin{split}
\frac{\partial \log L}{\partial b}=& \sum_{i=1}^{n}\frac{-a^2 c ((\beta -1) m+x_i)}{b (a+b) (a c+b (c+x_i)+m-x) (b (\beta -1) (c+x_i)+(\beta -1) m+x_i)}\\
 &+\sum_{i=1}^{n}\frac{a \left((\beta -1) b^2 (-x_i) (c+x_i)-2 b (c+x_i) ((\beta -1) m+x_i)+(x_i-m) ((\beta -1) m+x_i)\right)}{b (a+b) (a c+b (c+x_i)+m-x_i) (b (\beta -1) (c+x_i)+(\beta -1) m+x_i)}\\
 & -\sum_{i=1}^{n}\frac{b^2 \beta  x (c+x_i)}{b (a+b) (a c+b (c+x_i)+m-x_i) (b (\beta -1) (c+x_i)+(\beta -1) m+x_i)} 
\end{split}
\end{equation}
\begin{equation}
\begin{split}
\frac{\partial \log L}{\partial c}=&-\sum_{i=1}^{n}\frac{a (x_i (b (\beta -1)+1)+(\beta -1) m)+b \beta  x_i}{(a c+b (c+x_i)+m-x_i) (b (\beta -1) (c+x_i)+(\beta -1) m+x_i)}
\end{split}
\end{equation}
where $\psi(r)=\frac{d}{dr}\Gamma(r)$ is digamma function. As the above  equations  are not in closed form and hence cannot be solved explicitly. So we make use of a suitable  iterative technique to find the ML estimates  numerically like; by using \texttt{maxLik()} function in \texttt{R}. 
\section{Numerical Illustration}
\noindent  Table $1$   contains  automobile insurance data used by Michel Denuit($1997$). The proposed model i.e. Generalized Negative Binomial-Beta Exponential distribution(GNB-BE) is compared  with other models like Poisson  distribution(PD), Negative Binomial(NB) and Negative Binomial-Beta Exponential (NB-BE). The parameters are estimated by using the Maximum likelihood estimation.\\

\begin{table}[htbp]
	\centering
	\caption{ Distribution of Automobile insurance claim counts of Belgium 1958}
	\begin{tabular}{rrrrrrrr}
		\hline
		& Observed & \multicolumn{6}{c}{Expected} \\
		\hline
		Count & Frequency & PD      & NB     & NB-BE       & GNB-BE    \\
		0     & 7840   & 7635.62   & 7847.01   & 7845.72  & 7844.33  \\
		1     & 1317    & 1636.73   & 1288.36   & 1300.98 & 1300.24    \\
		2     & 239    & 175.419  & 256.533   & 243.204  & 243.706    \\
		3     & 42    & 12.5339  & 54.0665   & 52.7845   & 53.5029  \\
		4     & 14    & 0.671675  & 11.7097   & 13.0503  & 13.5028   \\
		5     & 4     & 0.028795  & 2.57699   & 3.6023  & 3.83641  \\
		6     & 4     & 0.001029  & 0.573128   & 1.09186   & 1.2057    \\
		7     & 1     & 0.000032   & 0.128418   & 0.358608   & 0.413233    \\
	
		Total & 9461   & 9461        & 9461       & 9461       &9461        \\
		&       &       &       &       &       &       &  \\ \hline
		& Estimated parameters &  ($\hat{\lambda}$) &  ($\hat{r},\hat{p}$)&  ($\hat{r},\hat{a},\hat{b},\hat{c}$) &   ($\hat{m},\hat{a},\hat{b},\hat{c},\hat{\beta}$)   &   \\
		&  &  0.214354 &  0.701512, 0.765955& 1.77119,1.96502  &  2.10804,1.9639  &   \\
		&  &   &  &  7.97405, 2.10205 &    4.1643,4.58003    &   \\
			&  &   &  &   &   1.09817   &   \\ \hline
		&       &       &       &       &       &       &  \\
		&       &       &       &       &       &       &  \\
		& log likelihood & -5490.78 & -5348.04 & -5343.80 & -5343.60  \\
		\hline
	\end{tabular}%
	\label{tab:addlabel}%
\end{table}
From the table it is evident as the log-likelihood of our proposed model is comparatively higher which means it outperforms other models.
\section{Conclusion}

In this article we introduce $\mathcal{GNB-BE}(m,\beta,a,b,c)$ distribution which is obtained by mixing  $ \mathcal{GNB}(m,\theta, \beta )$ distribution  and $\mathcal{BE}(a,b,c)$ distribution. We find that the  Negative binomial-beta exponential distribution, Generalized Waring distribution, Waring distribution and Yule distribution are all special cases of our proposed model. All the key moments of the  $\mathcal{GNB-BE}(m,\beta,a,b,c)$ distribution have been derived which includes factorial moments, mean and variance. Parameters have been estimated  Maximum Likelihood Estimation(MLE) method. The usage of the proposed model has been done by taking the real data into consideration


\begin{thebibliography}{}

\bibitem{} Abramowitz, M. and  Stegun, I. (1972). \textit{Handbook of Mathematical Function}, 2nd ed., Dover, New York.

\bibitem{} C. Pudprommarat, W. Bodhisuwan and P. Zeephongsekul(2012). A new mixed
negative binomial distribution,\textit{ J. Applied Sciences} 17 , 1853-1858.

\bibitem{} Deniz, E.,  Sarabia, J. and  Ojeda, E.(2008). Univariate and Multivariate Versions of the Negative Binomial-Inverse Gaussian Distributions with Applications, \textit{Insurance Mathematics and Economics} 42, 39-49.


\bibitem{} Denuit, M., (1997). A new distribution of Poisson-type for the numbers of claims. \textit{ Astin Bulletin} 27 (2), 229-242.

\bibitem{} E., Xekalaki (1983). A property of the Yule distribution and its applications. \textit{Commun. Stat. Theory Methods.} 12: 1181-1189.

\bibitem{} I., S., Gardshteyn and Ryzhik, I., M.(2007). \textit{Tables of Integrals, Series and Products} 7th Edn, Academic Press, New York.

\bibitem{} J., O., Irwin (1968). The generalized Waring distribution to accident theory. \textit{J. Roy. Stat. Soc. A,} 131: 205-225.

\bibitem{} J., O., Irwin (1975). The generalized Waring distribution Part. \textit{ I. J. Roy. Stat. Soc. A,} 138: 18-31. 

\bibitem{} Johnson, N. L., Kemp, A.W., and Kotz, S. (1992). \textit{Univariate discrete distributions}, New York:Wiley.

\bibitem{} Rahid,  A.,  Ahmad, Z. and Jan, T.(2014.) A Mixture of Generalized Negative Binomial Distribution with Generalized Exponential Distribution. \textit{Journal of Statistics Applications $\&$ Probability},3, 451-464.


\bibitem{} S., Nadarajah and Kotz, S.,(2006). The beta exponential distribution. \textit{Reliab. Eng. Syst. Saf.}, 91: 689-697.

\bibitem{} Wang, Z. (2011). One Mixed Negative Binomial Distribution with Application, \textit{Statistical Planning and Inference} 141 ,1153-1160.

\bibitem{}  Winkelmann, R.(2003). \textit{Econometric Analysis of Count Data}, 3rd Edn, Springer-Verlag, Berlin, Germany.



 
\end{thebibliography}
\end{document}